\documentclass[11pt,a4paper,reqno]{amsart}
\usepackage{amssymb,amsmath,amsthm}
\usepackage[utf8]{inputenc}
\usepackage[T1]{fontenc}
\usepackage{enumerate}
\usepackage[all]{xy}
\usepackage{fullpage}
\usepackage{comment}
\usepackage{array}
\usepackage{longtable}
\usepackage{stmaryrd}
\usepackage{mathrsfs} 
\usepackage{xcolor}
\usepackage{mathtools}

\def\F{\mathbb{F}}

\def\deg{\mathrm{deg}}

\def\dim{\mathrm{dim}}

\usepackage{hyperref}
\hypersetup{
    colorlinks=true,
    linkcolor=blue,
    filecolor=red,  
    citecolor=green,
    urlcolor=magenta,
    pdftitle={GON},
    pdfpagemode=FullScreen,
    }

\usepackage{cleveref}
\crefformat{section}{\S#2#1#3}
\crefformat{subsection}{\S#2#1#3}
\crefformat{subsubsection}{\S#2#1#3}
\usepackage{enumitem}
\usepackage{tikz}
\usepackage{mathdots}

\newtheorem{theorem}{Theorem}[section]
\newtheorem{lemma}[theorem]{Lemma}

\newtheorem{corollary}[theorem]{Corollary}
\newtheorem{proposition}[theorem]{Proposition}

\theoremstyle{definition}
\newtheorem{definition}[theorem]{Definition}

\theoremstyle{remark}
\newtheorem{remark}[theorem]{Remark}

\makeatletter
\newcommand{\subalign}[1]{%
  \vcenter{%
    \Let@ \restore@math@cr \default@tag
    \baselineskip\fontdimen10 \scriptfont\tw@
    \advance\baselineskip\fontdimen12 \scriptfont\tw@
    \lineskip\thr@@\fontdimen8 \scriptfont\thr@@
    \lineskiplimit\lineskip
    \ialign{\hfil$\m@th\scriptstyle##$&$\m@th\scriptstyle{}##$\hfil\crcr
      #1\crcr
    }%
  }%
}
\makeatother

\numberwithin{equation}{section}
\title{Geometric Properties of Periodic Lattices in Function Fields}
\author{Noy Soffer Aranov}
\email{noy.sofferaranov@tugraz.at}
\address{Graz University of Technology, Institute of Analysis and Number Theory, 8010 Graz, Austria}
\begin{document}
\maketitle

\begin{abstract}
    Periodic lattices are natural generalizations of lattices, which arise naturally in diophantine approximations with rationals of bounded denominators. In this paper, we prove analogues of classical theorems in geometry of numbers for periodic lattices in function fields. Moreover, we use special matrices to compute the covering and packing radii of special periodic lattices. 
\end{abstract}
\section{Introduction}

Let $d\geq 2$, let $S\subseteq \mathbb{R}^d$ be a Delone set, and let $\mathcal{C}$ be a $0$-symmetric convex body. Successive minima are geometric quantities associated to a Delone set. For $i=1,\dots,d$ and a $0$-symmetric convex body $\mathcal{C}$, we define the \emph{$i$-th successive minima} of $S$ with respect to $\mathcal{C}$ by 
\begin{equation}
\label{eqn:succMinDef}
    \lambda_{i,\mathcal{C}}(S)=\min\{r\geq 0:\dim(\operatorname{span}(S\cap B_{\Vert \cdot\Vert_{\mathcal{C}}}(0,r))\geq i\},
\end{equation}
where $B_{\Vert \cdot \Vert_{\mathcal{C}}}(0,r)=\{\mathbf{v}\in \mathbb{R}^d:\Vert \mathbf{v}\Vert_{\mathcal{C}}\leq r\}$ and $\Vert \mathbf{v} \Vert_{\mathcal{C}}=\inf\{r>0:\mathbf{v}\in r\mathcal{C}\}$. 

Lattices are a type of Delone set of the form $g\mathbb{Z}^d$, where $g\in \operatorname{GL}_d(\mathbb{R})$. They are very well studied, due to their rich algebraic and geometric structure. Moreover, the determinant of a lattice is proportional to the product of its successive minima.
\begin{theorem}[Minkowski's 2nd Theorem]
\label{thm:MinkReal}
For every lattice $\Lambda=g\mathbb{Z}^d$ and for every $0$-symmetric convex body $\mathcal{C}\subseteq \mathbb{R}^d$, we have
\begin{equation}
    \frac{2^n}{n!}\frac{\vert \det(g)\vert}{\operatorname{Vol}(\mathcal{C})}\leq \prod_{i=1}^d\lambda_{i,\mathcal{C}}(\Lambda)\leq 2^n\frac{\vert \det(g)\vert}{\operatorname{Vol}(\mathcal{C})},
\end{equation}
where $\operatorname{Vol}$ is the Lebesgue measure on $\mathbb{R}^d$.
\end{theorem}
Some other well studied geometric quantities associated with lattices and general Delone sets are packing and covering radii. 
\begin{definition}
\label{defs:Real}
Let $S$ be a Delone set, and let $\mathcal{C}\subseteq \mathbb{R}^d$ be a $0$-symmetric convex body.
    \begin{enumerate}
        \item The \emph{covering radius} of $S$ with respect to $\mathcal{C}$ is 
        $$\operatorname{CovRad}_{\mathcal{C}}(S)=\inf\{r>0|S+r\mathcal{C}=\mathbb{R}^d\}.$$
        \item The \emph{packing radius} of $S$ with respect to $\mathcal{C}$ is defined as $$\operatorname{PackRad}_{\mathcal{C}}(S)=\sup\left\{r>0:\min_{\mathbf{v}_1\neq\mathbf{v}_2\in S}\Vert \mathbf{v}_1-\mathbf{v}_2\Vert_{\mathcal{C}}>r\right\}.$$
        \item The set $S+\mathcal{C}$ is called a \emph{packing} if for every distinct $\mathbf{s}_1,\mathbf{s}_2$, we have $\mathbf{s}_1-\mathbf{s}_2\notin \mathcal{C}$. In this case, the \emph{upper density} of the packing is defined by
        $$\delta(S+\mathcal{C})=\limsup_{R\rightarrow \infty}\frac{\operatorname{Vol}\left((S+\mathcal{C})\cap B(0,R)\right)}{\operatorname{Vol}(B(0,R))}.$$
        We say that $S+\operatorname{PackRad}_{\mathcal{C}}(S)\mathcal{C}$ is the densest packing of $S$ with respect to $\mathcal{C}$.
        \item Given a convex body $\mathcal{C}$, define the packing density of $\mathcal{C}$ by
        $$\delta(\mathcal{C})=\sup\{\delta(S+\mathcal{C}):S\subseteq \mathbb{R}^d,S+\mathcal{C}\text{ is a packing}\}.$$
    \end{enumerate}
\end{definition}
Covering radii and packing radii of lattices are very well studied, due to their connection to cryptography, number theory, and diophantine approximations. For more information about covering radii of lattices in the real setting see \cite{Rog,RogCP,CS,RSW,ORW}. There has been significant progress about packings in recent years such as \cite{CS,RogCP,CE,V8,CKMRV}.

A natural generalization of lattices is periodic lattices, which are unions of finitely many translates of lattices. 
\begin{definition}
    Let $d\geq 2$, let $N\in \mathbb{N}$ and let $\Lambda \subseteq \mathbb{R}^d$ be a lattice. We say that a Delone set $S$ is a $(\Lambda,N)$ periodic lattice if:
    \begin{enumerate}
        \item $\Lambda+S=S$;
        \item Every fundamental cell of $\Lambda$ contains exactly $N$ points of $S$.
    \end{enumerate}
\end{definition}
Although abstract periodic lattices may be difficult to study, specific periodic lattices can provide useful information about diophantine approximations. For example, Aliev and Henk \cite{AH} studied the successive minima of certain periodic lattices to study approximations of vectors by rational vectors with bounded denominators. For a lattice $\Lambda=g\mathbb{Z}^n$, vectors $\alpha,\mathbf{v}\in \mathbb{R}^n$, a number $N\in \mathbb{N}$, and a convex body $\mathcal{C}\subseteq \mathcal{R}^n$, define
$$\beta(\Lambda,\alpha,\mathcal{C},N)=\left\{\Vert Q\alpha-\mathbf{v}\Vert_{\mathcal{C}}:\mathbf{v}\in \Lambda,Q\in \{0,1,\dots,N\}\right\}.$$
Assume that for every $Q\in \{1,\dots,N\}$, we have $Q\alpha\notin \Lambda$. In order to study the set $\beta(\Lambda,\alpha,\mathcal{C},N)$, we define the $i$-th successive minima of the tuple $(\Lambda,\alpha,\mathcal{C},N)$, for $i=1,\dots,d$ by
\begin{equation}
    \tilde{\lambda}_i(\Lambda,\alpha,\mathcal{C},N)=\min\left\{r>0:\begin{matrix}\text{there exist }i\text{ linearly independent vectors }\\(Q_j,\mathbf{v}_j)\in \{0,1,\dots,N\}\times \Lambda,\\ \text{ such that }\Vert Q_j\alpha-\mathbf{v}_j\Vert_{\mathcal{C}}\leq r\end{matrix}\right\}.
\end{equation}
Aliev and Henk \cite{AH} connected between the successive minima $\tilde{\lambda}_i$ and successive minima of certain periodic lattices. For a lattice $\Lambda$, a vector $\alpha\in \mathbb{R}^d$, and $N\in \mathbb{N}$, define
\begin{equation}
    \Lambda(\alpha,N)=\bigcup_{Q\in \{0,1,\dots,N\}}(Q\alpha+\Lambda).
\end{equation}
It is clear that
\begin{equation}
    \tilde{\lambda}_i(\Lambda,\alpha, \mathcal{C},N)\leq \lambda_{i,\mathcal{C}}(\Lambda(\alpha,N)). 
\end{equation}
Therefore, to study the set $\beta(\Lambda,\alpha,\mathcal{C},N)$, Aliev and Henk bounded $\lambda_{i,\mathcal{C}}(\Lambda(\alpha,N))$.
\begin{theorem}{\cite[Theorem 1.1]{AH}}
    Let $\Lambda$ be a lattice, let $\alpha\in\mathbb{R}^d$ be a vector, let $\mathcal{C}$ be a convex body and let $N\in \mathbb{N}$. Then, 
    \begin{equation}
        \lambda_{1,\mathcal{C}}(\Lambda(\alpha,N))^d\leq \delta(\mathcal{C})2^d\frac{\det(\Lambda)}{N+1}.
    \end{equation}
    Moreover, 
    \begin{equation}
        \frac{2^d}{d!}\det(\Lambda)\gamma(\alpha,\Lambda,N,d)\leq \prod_{i=1}^d\lambda_{i,\mathcal{C}}(\Lambda(\alpha,N))\leq \delta(\mathcal{C})2^d\frac{\det(\Lambda)}{N+1},
    \end{equation}
    where $\gamma(\alpha,\Lambda,N,d)$ is an explicit geometric constant depending on $\alpha,\Lambda,N,$ and $d$.
\end{theorem}
In this paper, we study geometric properties of $\Lambda(\alpha,N)$ and general periodic lattices in the function field setting, such as successive minima, covering radius, packing radius, and packing density. Several of these results are generalizations of the results of \cite{Mah,BK,A} or function field analogues of the results of \cite{AH}. In recent years, there have been several advancements in geometry of numbers in function fields such as \cite{KST,GLS24,BKL,AB24b,AK,Bang,A,BK}. 
\begin{remark}
    Geometry of periodic lattices can also be studied for periodic lattices defined over general global fields. In such fields, one loses the exact bounds, since Minkowski's 2nd Theorem (Theorem \ref{thm:Mink2nd}) does not yield an equality for global fields. Hence, in this paper, we focus only on the function field setting.
\end{remark}
\subsection{The Function Field Setting}
Let $p$ be a prime, and let $q$ be a power of $p$. Let $\mathcal{R}=\mathbb{F}_q[x]$, and for $N\in \mathbb{N}$, let $\mathcal{R}_{\leq N}=\{Q\in \mathbb{F}_q[x]:\deg(Q)\leq N\}$. For a rational function $\frac{f}{g}\in \mathbb{F}_q(x)$, define the absolute value of $\frac{f}{g}$ by $\left|\frac{f}{g}\right|=q^{\deg(f)-\deg(g)}$. Let $\mathcal{K}_{\infty}=\mathbb{F}_q((x^{-1}))$ be the completion of $\mathcal{K}=\mathbb{F}_q(x)$ with respect to the absolute value $\vert \cdot \vert$. For $r>0$, denote $B(0,r)=\{\alpha\in \mathcal{K}_{\infty}:\vert \alpha\vert\leq r\}$. Let $\mathcal{O}=\mathbb{F}_q[[x^{-1}]]=B(0,1)$ be the maximal compact order. Define the maximal ideal of $\mathcal{O}$ by $\mathfrak{m}=x^{-1}\mathcal{O}=B(0,q^{-1}).$

Let $d\geq 2$. Define the Haar measure on $\mathcal{K}_{\infty}^d$ to be the unique translation invariant measure such that $m(\mathcal{O}^d)=1$. As seen in \cite{Mah}, a convex body in $\mathcal{K}_{\infty}^d$ is a set of the form $g\mathcal{O}^d$, where $g\in \operatorname{GL}_d(\mathcal{K}_{\infty})$. Then, the volume of the convex body $\mathcal{C}=g\mathcal{O}^d$ is $m(g\mathcal{O}^d)=\vert \det(g)\vert$. A convex body $\mathcal{C}$ induces a norm on $\mathcal{K}_{\infty}^d$ by $$\Vert \mathbf{v}\Vert_{\mathcal{C}}=\min\{r>0:\mathbf{v}\in q^{\log_qr}\mathcal{C}\}.$$ In particular, the convex body $\mathcal{O}^d$ induces the supremum norm $\Vert \mathbf{v}\Vert=\max_{i=1,\dots,d}\vert v_i\vert$. Both the absolute value $\vert\cdot \vert$ and the supremum norm $\Vert \cdot \Vert$ satisfy the ultrametric inequality.
\begin{lemma}{Ultrametric Inequality}
\label{lem:UM}
    \begin{enumerate}
            \item For all $\alpha,\beta\in \mathbb{F}_q((x^{-1}))$, we have $\vert \alpha+\beta\vert\leq \max\{\vert \alpha\vert,\vert \beta\vert\}$.
            \item For all $\mathbf{u},\mathbf{v}\in \mathbb{F}_q((x^{-1}))^n$, we have $\Vert\mathbf{u}+\mathbf{v}\Vert\leq \max\{\Vert \mathbf{u}\Vert,\Vert\mathbf{v}\Vert\}$. 
        \end{enumerate}        
\end{lemma}
\subsection{Lattices and Periodic Lattices in $\mathcal{K}_{\infty}^d$}
\begin{definition}
A lattice in $\mathcal{K}_{\infty}^d$ is a set of the form $g\mathcal{R}^d$, where $g\in \operatorname{GL}_d(\mathcal{K}_{\infty})$. A convenient fundamental domain for a lattice $\Lambda=g\mathcal{R}^d$ is $\mathcal{D}_{\Lambda}=g\mathfrak{m}^d$. Motivated by this, we define $\det(\Lambda)=\vert \det(g)\vert$, so that $\det(\Lambda)=q^dm(\mathcal{D}_{\Lambda})$. 
\end{definition}
Lattices are a special type of Delone set, which are very well studied. Throughout this paper, we discuss the following quantities relating to Delone sets in $\mathcal{K}_{\infty}^d$, which are function field analogues of Definition \ref{defs:Real}. 
\begin{definition}
\label{def:GeomVals}
    Let $S\subseteq \mathcal{K}_{\infty}^d$ be a Delone set, and let $\mathcal{C}$ be a convex body.
    \begin{enumerate}
        \item The \emph{covering radius} of $S$ with respect to $\mathcal{C}$ is
        $$
\operatorname{CovRad}_{\mathcal{C}}(S)=\min\{ r > 0 \mid  S+x^{\log_qr}\mathcal{C}=\mathcal{K}_{\infty}^d \}=\sup_{\mathbf{v}\in\mathcal{K}_{\infty}^d}\inf_{\mathbf{u}\in S}\Vert \mathbf{v}-\mathbf{u}\Vert_{\mathcal{C}}.
        $$
        When $\mathcal{C}=\mathcal{O}^d$, we write $\operatorname{CovRad}:=\operatorname{CovRad}_{\mathcal{O}^d}$. 
        \item The \emph{packing radius} of $S$ with respect to $\mathcal{C}$ is
        
        $$\operatorname{PackRad}_{\mathcal{C}}(S)=\max\left\lbrace r\geq 0 \mid \min_{{\bf s}_1 \neq {\bf s}_2 \in S}\Vert {\bf s}_1-{\bf s}_2 \Vert_{\mathcal{C}}>r \right\rbrace$$
        $$=\max\left\{r\geq 0:\forall \mathbf{s}_1,\mathbf{s}_2\in S:\mathbf{s}_1\neq \mathbf{s}_2,(\mathcal{C}+\mathbf{s}_1)\cap (\mathcal{C}+\mathbf{s}_2)=\emptyset\right\}.
        $$ 
        Similarly, when, $\mathcal{C}=\mathcal{O}^d$, we write $\operatorname{PackRad}:=\operatorname{PackRad}_{\mathcal{O}^d}$.
        \item Let $\mathcal{C}$ be a convex body. The set $S+\mathcal{C}$ is called a \emph{packing} if, for every pair of distinct points $\mathbf{s}_1, \mathbf{s}_2 \in S$, we have
        $
        \mathbf{s}_1 - \mathbf{s}_2 \notin \mathcal{C}
        $. 
        In this case, the \emph{upper density} of the packing is defined by 
        $$\delta(S+\mathcal{C})=\limsup_{R\rightarrow \infty}\frac{m((S+\mathcal{C})\cap B(0,R))}{m(B(0,R))}.$$ 
        The packing $S+x^{\log_q\operatorname{PackRad}_{\mathcal{C}}(S)}\mathcal{C}$ is called the densest packing of $S$ with respect to $\mathcal{C}$.
        \item For $i=1,\dots,d$, the \emph{$i$-th successive minima} of $S$ with respect to $\mathcal{C}$ is defined by
        $$\lambda_{i,\mathcal{C}}(S)=\min\{r>0 \mid \dim(\operatorname{span}_{\mathcal{K}_\infty}(S\cap B_{\Vert\cdot\Vert_{\mathcal{C}}}(0,r)))=i\}.$$
        When $\mathcal{C}=\mathcal{O}^d$, we write $\lambda_i:=\lambda_{i,\mathcal{O}^d}$.
        \item If $\mathbf{v}_1,\dots,\mathbf{v}_d$ are a set of independent vectors such that $\lambda_{i,\mathcal{C}}(S)=\Vert \mathbf{v}_i\Vert_{\mathcal{C}}$, then, we say that $\mathbf{v}_1,\dots,\mathbf{v}_d$ are a set of successive minima for $S$.
    \end{enumerate}
\end{definition}
Successive minima, covering radii, and packing radii are well-studied geometric quantities for lattices in $\mathcal{K}_{\infty}^d$. In the function field setting, there is an analogue of Theorem \ref{thm:MinkReal}.
\begin{theorem}[Minkowski's 2nd theorem in function fields \cite{Mah}]
\label{thm:Mink2nd}
    Let $\Lambda\subseteq \mathcal{K}_{\infty}^d$ be a lattice and let $\mathcal{C}$ be a convex body. Then, 
    \begin{equation}
        \prod_{i=1}^d\lambda_{i,\mathcal{C}}(\Lambda)=\frac{\det(\Lambda)}{m(\mathcal{C})}.
    \end{equation}
\end{theorem}
Furthermore, every lattice $\Lambda$ has an $\mathcal{R}$ basis of successive minima.
\begin{lemma}{\cite[Lemma 1.18]{A}}
    Let $\Lambda$ be a lattice, and let $\mathbf{v}^{(1)},\dots,\mathbf{v}^{(d)}\in \Lambda$ be linearly independent vectors satisfying $\Vert \mathbf{v}^{(i)}\Vert=\lambda_i(\Lambda)$. Then, $\Lambda=\operatorname{span}_{\mathcal{R}}\{\mathbf{v}^{(1)},\dots,\mathbf{v}^{(d)}\}$.
\end{lemma}
Moreover, the covering radius of a lattice has a closed formula.
\begin{theorem}{\cite[Corollary 1.26]{A}}
    Let $\Lambda\subseteq \mathcal{K}_{\infty}^d$ be a lattice and let $\mathcal{C}$ be a convex body. Then, $\operatorname{CovRad}_{\mathcal{C}}(\Lambda)=q^{-1}\lambda_{d,\mathcal{C}}(\Lambda)$
    
\end{theorem}
In this paper, we study geometry of numbers for specific period lattices, arising from a classic question is diophantine approximations, which discusses approximation of a vector by rational vectors with bounded denominators. Toward this end, we make the following definition.
\begin{definition}
    Let $N\in \mathbb{N}$, let $\Lambda$ be a lattice, and let $\mathbf{v}\in \mathcal{K}_{\infty}^n$ be a vector. We say that $\mathbf{v}$ is $N$-rational with respect to $\Lambda$ if there exists $Q\in \mathcal{R}_{\leq N}\setminus \{0\}$ such that $Q\mathbf{v}\in \Lambda$. Otherwise, we say that $\mathbf{v}$ is $N$-irrational with respect to $\Lambda$.
\end{definition}
In particular if $\mathcal{C}\subseteq \mathcal{K}_{\infty}^n$ is a convex body, $N\in \mathbb{N}$, $\Lambda$ is a lattice, and $\alpha,\mathbf{v}\in \mathcal{K}_{\infty}^d$ is $N$-irrational with respect to $\Lambda$, one can study
$$\lambda_{i,\mathcal{C}}(\alpha,\mathbf{v},\Lambda,N)=\min\left\{r>0:\begin{array}{ccc}
     \exists \mathbf{P}_1,\dots,
     \mathbf{P}_{i}\in \Lambda,\,Q_1,\dots,Q_{i}\in \mathcal{R}_{\leq N},\\
     \dim\operatorname{span}\left\{(\mathbf{P}_j,Q_j):j=1,\dots,i\right\}\geq i,\\
     Q_i\alpha-\mathbf{v}-\mathbf{P}\in x^{\log_qr}\mathcal{C}.
\end{array}\right\}.$$
A useful framework to study homogeneous (when $\mathbf{v}=0$) and inhomogeneous approximation with bounded denominators is through geometric properties of special periodic lattices defined by
$$\Lambda(\alpha,q^N)=\bigcup_{\deg(Q)\leq N}(Q\alpha+\Lambda).$$
Then, $\lambda_{1,\mathcal{C}}(\alpha,0,\Lambda,N)=\lambda_{1,\mathcal{C}}(\Lambda(\alpha,q^N))$ and $\lambda_{i,\mathcal{C}}\left(\alpha,0,\Lambda,N\right)\leq \lambda_{i,\mathcal{C}}\left(\Lambda(\alpha,q^N)\right)$ for every $i=2,\dots,d$, where $\lambda_{i,\mathcal{C}}(\Lambda(\alpha,q^N)))$ are defined in Definition \ref{def:GeomVals}. In this paper, we bound and compute the successive minima and other geometric quantities associated to $\Lambda(\alpha,q^N)$. We first define general periodic lattices, which are generalizations of sets of the form $\Lambda(\alpha,q^N)$.
\begin{definition}\label{def:periolatt}
    Let $\Lambda \in \mathcal{L}_d$ be a lattice and let $N\in \mathbb{N}$. An $\F_q$-subspace $S \subset \mathcal{K}_\infty^d$ is called a \emph{$(\Lambda,q^N)$-periodic lattice} if:
    \begin{enumerate}
        \item \label{LambdaInv}$\Lambda+S=S$;
        \item \label{D_LambdaCapS-Size} $\#\mathcal{D}_{\Lambda}\cap S=q^N$.
    \end{enumerate}
\end{definition}
In particular, every $(\Lambda,q^N)$-periodic lattice $S$ contains $\Lambda$. For example, $\Lambda$ is a $(\Lambda,q^0)$-periodic lattice. Moreover, if $\alpha$ is irrational, then $\Lambda(\alpha,q^N)$ is a $(\Lambda,q^{N+1})$-periodic lattice (see Lemma \ref{lem:intersectdlamb}). 
\begin{remark}
    When $S=\Lambda(\alpha,q^N)$, we have $\operatorname{CovRad}_{\mathcal{C}}(S)=\sup_{\mathbf{v}\in \mathcal{K}_{\infty}^n}\lambda_{1,\mathcal{C}}(\alpha,\mathbf{v},\Lambda,N)$.
\end{remark}
\subsection{Main Results}
We prove an analogue of Minkowski's convex body theorem \cite{Mah},\cite[Theorem 5.5]{Cla} for periodic lattices. 
\begin{proposition}[Minkowski's convex body theorem for periodic lattices]
\label{prop:Mink}
    Let $S$ be a $(\Lambda,q^N)$ periodic lattice, and let $\mathcal{C}$ be a convex body satisfying $m(\mathcal{C}+\mathcal{D}_{\Lambda}\cap S)>\frac{\det(\Lambda)}{q^{N+d}}$. Then, $\mathcal{C}$ contains a non-zero point of $S$. 
\end{proposition}
In order to prove Proposition \ref{prop:Mink}, we prove a lattice point counting claim, which is the periodic lattice analogue of \cite[Lemma 6.2]{BK}. Furthermore, one can view this as a function field analogue of a conjecture of Betke, Henk, and Willis \cite[Conjecture 2.1]{BHW} (see also \cite[Conjecture 2.1]{Toi})
\begin{proposition}
\label{prop:BKPerLatt}
    Let $\mathcal{C}$ be a convex body, let $\Lambda$ be a lattice, let $N\in \mathbb{N}$, and let $S$ be a $(\Lambda,q^N)$ periodic lattice. Then, 
    \begin{equation}
        \left|\mathcal{C}\cap S\right|=\left|\mathcal{C}\cap \mathcal{D}_{\Lambda}\cap S\right|\prod_{i=1}^d\left\lceil \frac{q}{\lambda_{i,\mathcal{C}}(\Lambda)}\right\rceil.
    \end{equation}
    In particular, if $\mathcal{D}_{\Lambda}\subseteq \mathcal{C}$, then, 
        \begin{equation}
            \vert \mathcal{C}\cap S\vert=q^N\prod_{i=1}^d\left\lceil \frac{q}{\lambda_{i,\mathcal{C}}(\Lambda)}\right\rceil.
        \end{equation}
\end{proposition}
Moreover, we compute the packing radius and packing density of the densest packing.
\begin{lemma}
\label{lem:PackRadComp}
    For every periodic lattice $S\subseteq \mathcal{K}_{\infty}^d$ and for every convex body $\mathcal{C}$, we have $\operatorname{PackRad}_{\mathcal{C}}(S)=q^{-1}\lambda_{1,\mathcal{C}}(S)$.
\end{lemma}

\begin{proposition}
\label{lem:PackDens}
    We have  
    $$\delta(x^{-1+\log_q\lambda_{1,\mathcal{C}}(S)}\mathcal{C}+S)=q^N\lambda_{1,\mathcal{C}}(S)^d\frac{m(\mathcal{C})}{\det(\Lambda)}.$$ 
    When $S=\Lambda$, we have $\delta(x^{-1+\log_q\lambda_{1,\mathcal{C}}(\Lambda)}\mathcal{C}+\Lambda)=\frac{\lambda_{1,\mathcal{C}}(\Lambda)^d}{\det(\Lambda)}$.
\end{proposition}
As a corollary, we bound $\lambda_{1,\mathcal{C}}(S)$, and in addition, we bound the product of the successive minima. 
\begin{theorem}
\label{thm:SuccMinBNDS}
    Let $\Lambda\subseteq \mathcal{K}_{\infty}^d$ be a lattice, let $N\in \mathbb{N}$, let $\mathcal{C}$ be a convex body, and let $S$ be a $(\Lambda,q^N)$ periodic lattice. Then, 
    \begin{enumerate}
        \item \label{thm:lambda_1BND} $\lambda_{1,\mathcal{C}}(S)^d\leq \frac{\det(\Lambda)}{q^Nm(\mathcal{C})}$;
        \item \label{thm:prodBNDPer} $\prod_{i=1}^d\lambda_{i,\mathcal{C}}(S)\leq \frac{\det(\Lambda)}{q^Nm(\mathcal{C})}$.
    \end{enumerate}
\end{theorem}
\begin{proof}[Proof of Theorem \ref{thm:SuccMinBNDS}(\ref{thm:lambda_1BND})]
By Proposition \ref{lem:PackDens}, we have $q^N\lambda_{1,\mathcal{C}}(S)^d\frac{m(\mathcal{C})}{\det(\Lambda)}\leq 1$, so that Theorem  \ref{thm:SuccMinBNDS}(\ref{thm:lambda_1BND}) follows. 
\end{proof}
When $S=\Lambda(\alpha,q^N)$, we obtain a lower bound on the product of the successive minima as well, which can be viewed as an analogue of \cite[Theorem 1.1]{AH}.
\begin{theorem}
\label{thm:SuccMinProd}
    Let $\Lambda\subseteq \mathcal{K}_{\infty}^d$ be a lattice, let $\mathbf{v}^{(1)},\dots,\mathbf{v}^{(d)}$ be a basis of successive minima for $\Lambda$, let $\mathcal{C}$ be a convex body, and let $N\in \mathbb{N}$. Let $\alpha=\sum_{i=1}^d\alpha_i\mathbf{v}^{(i)}$ be such that $\mathcal{R}_{\leq N}\alpha_i\cap\mathcal{R}=\{0\}$ for every $i=1,\dots ,d$. Then, 
    \begin{equation}
    \label{eqn:succMinIneq}
        d(\Lambda(\alpha,q^N))\frac{\det(\Lambda)}{m(\mathcal{C})}\leq \prod_{i=1}^d\lambda_{i,\mathcal{C}}\left(\Lambda(\alpha,q^N)\right)\leq \frac{\det(\Lambda)}{q^{N+1}m(\mathcal{C})},
    \end{equation}
    where 
    $$d(\Lambda(\alpha,q^N))=\min_{k=1,\dots,d}\min_{\{j_1,\dots,j_k\}\subseteq \{1,\dots,d\}}\left\{\left|\det\left(\langle Q_i\alpha_{j_{\ell}}\rangle\right)_{1\leq i,\ell\leq k,}\right|>0:\deg(Q_i)\leq N,\forall i=k+1,\dots,d\right\}.$$
    As a consequence, there exist $k=1,\dots,d$ and polynomials $Q_{k+1},\dots,Q_d\in \mathcal{R}_{\leq N}$ such that
    \begin{equation}
        0\neq \left|\det\begin{pmatrix}
        \langle Q_{1}\alpha_{\ell_1}\rangle&\dots&\langle Q_k\alpha_{\ell_1}\rangle\\
        \vdots&\dots&\vdots\\
        \langle Q_{1}\alpha_{\ell_k}\rangle&\dots&\langle Q_k\alpha_{\ell_k}\rangle
    \end{pmatrix}\right|\leq \frac{1}{q^{N+1}}.
    \end{equation}
\end{theorem}
In addition, we generalize \cite[Corollary 1.26]{A} to periodic lattices. To explain our result, we need some terminology about Hankel matrices. For $\alpha=\sum_{n=1}^{\infty}\alpha_nx^{-n}\in \mathfrak{m}$, we define the Hankel matrix of $\alpha$ of order $m\times n$ by 
\begin{equation}
    \Delta_{\alpha}(m,n)=\begin{pmatrix}
        \alpha_1&\alpha_2&\dots&\alpha_n\\
        \alpha_2&\alpha_3&\dots&\alpha_{n+1}\\
        \vdots&\dots&\ddots&\vdots\\
        \alpha_m&\dots&\dots&\alpha_{m+n-1}
    \end{pmatrix}.
\end{equation}
If $m\leq 0$ or $n\leq 0$, we let $\Delta_{\alpha}(m,n)$ be an empty matrix. These matrices play an important role in diophantine approximations over function fields \cite{dMT,ALN,GaGh17,A,AMin,R,GR,ANS,AR}, which we explain later. To compute the covering radius of a periodic lattice, it is useful to concatonate Hankel matrices. For $d\in \mathbb{N}$, $(\ell_1,\dots,\ell_d)\in \mathbb{Z}^d$, $(\alpha_1,\dots,\alpha_d)\in \mathfrak{m}^d$, and $N\in \mathbb{N}$, denote
\begin{equation}
    \Delta_{\alpha_1,\dots,\alpha_d}(\ell_1,\dots,\ell_d;N)=\begin{pmatrix}
        \Delta_{\alpha_1}(\ell_1,N)\\
        \vdots\\
        \Delta_{\alpha_d}(\ell_d,N)
    \end{pmatrix}.
\end{equation}
\begin{proposition}
\label{prop:CovRad}
    Let $\Lambda$ be a lattice whose successive minima are $\mathbf{v}^{(1)},\dots,\mathbf{v}^{(d)}$. Let $N\in \mathbb{N}$, let $\mathcal{C}$ be a convex body, and let $\alpha=\sum_{i=1}^d\alpha_i\mathbf{v}^{(i)}\in \mathfrak{m}$. Then, 
   \begin{equation}
        \operatorname{CovRad}_{\mathcal{C}}(\Lambda(\alpha,q^N))=q^{-\left(1+\gamma(\Lambda(\alpha,q^N))\right)},
    \end{equation}
    where $\gamma(\Lambda(\alpha,q^N))$ is the smallest value of $\ell\geq 0$, such that
    \begin{equation}
        \operatorname{rank}(\Delta_{\alpha_1,\dots,\alpha_d}(\ell+\log_q\lambda_{1,\mathcal{C}}(\Lambda),\dots,\ell+\log_q\lambda_{d,\mathcal{C}}(\Lambda);N+1))=\sum_{i=1}^d\max\{\ell+\log_q\lambda_{i,\mathcal{C}}(\Lambda),0\}.
    \end{equation}
\end{proposition}
\begin{remark}
    Note that \cite[Corollary 1.26]{A} can be obtained as a corollary of Proposition \ref{prop:CovRad}, when $\alpha=0$ and $N=0$. This stems from the fact that
    $$\sum_{i=1}^d\max\{\ell+\log_q\lambda_{d,\mathcal{C}}(\Lambda),0\}>0$$ if and only if $\ell>-\log_q\lambda_{i,\mathcal{C}}(\Lambda)$. The proof of Proposition \ref{prop:CovRad} uses systems of linear equations over $\mathbb{F}_q$, whereas the proof of \cite[Corollary 1.26]{A} uses reduction theory. 
\end{remark}
\begin{corollary}
    For every lattice $\Lambda\subseteq \mathcal{K}_{\infty}^d$ and for every $\alpha\in\mathcal{K}_{\infty}^d$, we have
    \begin{equation}
        q^{-\left(1+\max_{i=1,\dots,d}\frac{N+1-\sum_{j=d-i+1}^d\log_q\lambda_j(\Lambda)}{i}\right)}\leq \operatorname{CovRad}_{\mathcal{C}}(\Lambda(\alpha,q^N))\leq q^{-1}\lambda_{d,\mathcal{C}}(\Lambda)
    \end{equation}
\end{corollary}
\section{Proofs}
\subsection{Packing Radii}
\begin{proof}[Proof of Lemma \ref{lem:PackRadComp}]
Since $S$ is an $\mathbb{F}_q$ subspace and $\lambda_{1,\mathcal{C}}(S)=\min_{{\bf s}\in S \setminus \{{\bf 0}\}} \Vert{\bf s} \Vert_{\mathcal{C}}$, then
$$
\operatorname{PackRad}_{\mathcal{C}}(S)=\max\left\lbrace r \geq 0 \mid \min_{{\bf s}\in S \setminus \{{\bf 0}\}} \Vert{\bf s} \Vert_{\mathcal{C}} >r\right\rbrace=q^{-1}\lambda_{1,\mathcal{C}}(S).$$
\end{proof}
\begin{lemma}\label{lem:intersectdlamb}
   Let $\Lambda$ be a lattice, let $N\in \mathbb{N}$, and let $\alpha$ be a vector which is $N$-irrational with respect to $\Lambda$. Then, we have $\#\mathcal{D}_{\Lambda}\cap \Lambda(\alpha,q^N)=q^{N+1}$. As a consequence, $\Lambda(\alpha,q^N)$ is a $(\Lambda,q^{N+1})$ periodic lattice. 
\end{lemma}
\begin{proof}
    For every $Q\in \mathcal{R}_{\leq N}$, there exist unique $\vert \alpha_i\vert\leq q^{-1}$ ($i=1,\dots,d$) and $a_i\in \mathcal{R}$ ($i=1,\dots,d$), such that $Q\alpha=\sum_{i=1}^d(\alpha_i+a_i)\mathbf{v}_i$. Therefore, $Q\alpha-\sum_{i=1}^da_i\mathbf{v}_i=\sum_{i=1}^d\alpha_i\mathbf{v}_i\in \mathcal{D}_{\Lambda}$ is the unique vector of the form $Q\alpha-\mathbf{w}\in \mathcal{D}_{\Lambda}$, where $\mathbf{w}\in \Lambda$. Consequently, $$\#\mathcal{D}_{\Lambda}\cap \Lambda(\alpha,q^N)\leq \#\mathcal{R}_{\leq N}=q^{N+1}.$$
    Let $Q\alpha=\sum_{i=1}^d(a_i+\alpha_i)\mathbf{v}_i$ and $P\alpha=\sum_{i=1}^d(b_i+\beta_i)\mathbf{v}_i$, where $P,Q\in \mathcal{R}_{\leq N}$, $a_i,b_i\in \mathcal{R}$, and $\alpha_i,\beta_i\in \mathfrak{m}$. If $\sum_{i=1}^d\alpha_i\mathbf{v}_i=\sum_{i=1}^d\beta_i\mathbf{v}_i$, then, $(Q-P)\alpha=\sum_{i=1}^d(a_i-b_i)\mathbf{v}_i\in \Lambda$. Thus $\alpha$ is $N$-rational with respect to $\Lambda$, which contradicts the assumption that $\alpha$ is $N$-irrational with respect to $\Lambda$. This contradiction proves that the vectors $\sum_{i=1}^d\alpha_i\mathbf{v}_i$ are distinct, so that $\#\mathcal{D}_{\Lambda}\cap \Lambda(\alpha,q^N)=q^{N+1}$. Hence, $\Lambda(\alpha,q^N)$ is a $(\Lambda,q^{N+1})$ periodic lattice.  
\end{proof}
\begin{lemma}
\label{lem:PerLattCapBallCnt}
Let $\Lambda$ be a lattice, let $N\in \mathbb{N}$, and let $S$ be a $(\Lambda,q^N)$-periodic lattice. For every integer $R\geq 0$, we have
\begin{equation}
    \left|B(0,q^R)\cap S\right|=\left|\mathcal{D}_{\Lambda}\cap S\cap B(0,q^R)\right|\cdot\prod_{i=1}^d\left\lceil\frac{q^{R+1}}{\lambda_i(\Lambda)}\right\rceil.
\end{equation}
In particular, if $R$ is large enough so that $S\cap \mathcal{D}_{\Lambda}\subseteq B(0,q^R)$, 
\begin{equation}
    \left|B(0,q^R)\cap S\right|=q^N\prod_{i=1}^d\left\lceil\frac{q^{R+1}}{\lambda_i(\Lambda)}\right\rceil.
\end{equation}
\end{lemma}
Proposition \ref{prop:BKPerLatt} is an immediate corollary of Lemma \ref{lem:intersectdlamb} and Lemma \ref{lem:PerLattCapBallCnt}.
\begin{proof}[Proof of Proposition \ref{prop:BKPerLatt}]
    If $\mathcal{C}=h\mathcal{O}^d$, then, $h^{-1}S$ is a $(h^{-1}\Lambda,q^N)$ periodic lattice. Thus, by Lemma \ref{lem:PerLattCapBallCnt},
    \begin{align*}
        \left|\mathcal{C}\cap S\right|&=\left|h\mathcal{O}^d\cap S\right|\\
        &=\left|\mathcal{O}^d\cap h^{-1}S\right|\\
        &=\left|\mathcal{O}^d\cap h^{-1}\mathcal{D}_{\Lambda}\cap h^{-1}S\vert \right|\prod_{i=1}^d\left\lceil\frac{q}{\lambda_i(h^{-1}\Lambda)}\right\rceil\\
        &=\vert \mathcal{C}\cap \mathcal{D}_{\Lambda}\cap S\vert\prod_{i=1}^d\left\lceil \frac{q}{\lambda_{i,\mathcal{C}}(\Lambda)}\right\rceil.
    \end{align*}
    In particular, if $\mathcal{D}_{\Lambda}\subseteq \mathcal{C}$, then, $\left|\mathcal{C}\cap \mathcal{D}_{\Lambda}\cap S\right|=\vert S\cap \mathcal{D}_{\Lambda}\vert=q^N$, so that
    \begin{equation}
        \vert \mathcal{C}\cap \mathcal{S}\vert=q^N\prod_{i=1}^d\left\lceil \frac{q}{\lambda_{i,\mathcal{C}}(\Lambda)}\right\rceil.
    \end{equation}
\end{proof}
\begin{proof}[Proof of Lemma \ref{lem:PerLattCapBallCnt}]
    By \eqref{LambdaInv} and the fact that $\mathcal{K}_{\infty}^d=\mathcal{D}_\Lambda+\Lambda$, we have $S=\mathcal{D}_{\Lambda}\cap S+\Lambda$. Thus, elements $S\cap B(0,q^R)$ are of the form $\mathbf{v}+\mathbf{u}$, where 
    \begin{enumerate}
        \item $\mathbf{u}\in \mathcal{D}_{\Lambda}\cap B(0,q^R)\cap S$,
        \item $\mathbf{v}\in \Lambda$,
        \item \label{Cond:|u+v|<q^R}and $\Vert \mathbf{u}+\mathbf{v}\Vert\leq q^R$.
    \end{enumerate}
    By the Lemma \ref{lem:UM} -- observing that either ${\bf v}={\bf 0}$ or $\Vert \mathbf{v}\Vert>\Vert \mathbf{u}\Vert$ -- we may replace condition \eqref{Cond:|u+v|<q^R} with $\Vert \mathbf{v}\Vert\leq q^R$. Therefore, we have a direct sum 
    \begin{equation*}
    \label{eqn:B(0,q^R)capLambda(alpha,q^N)}
        B(0,q^R)\cap S=B(0,q^R)\cap \Lambda \oplus \mathcal{D}_{\Lambda}\cap S\cap B(0,q^R),
    \end{equation*}
    where the symbol $\oplus$ denotes a direct sum of the sets. Consequently, by \cite[Lemma 6.2]{BK},
\begin{align*}\label{eqn:PerLatcapBallPtCnt}
\left|B(0,q^R)\cap S\right|&=\left|\mathcal{D}_{\Lambda}\cap S\cap B(0,q^R)\right|\cdot \left|B(0,q^R)\cap \Lambda\right|\\
&=\left|\mathcal{D}_{\Lambda}\cap S\cap B(0,q^R)\right|\cdot \left|B(0,1)\cap x^{-R}\Lambda\right|\\
&=\left| \mathcal{D}_{\Lambda}\cap S\cap B(0,q^R)\right|\prod_{i=1}^d\left\lceil\frac{q}{\lambda_i(x^{-R}\Lambda)}\right\rceil\\
&=\left|\mathcal{D}_{\Lambda}\cap S\cap B(0,q^R)\right|\prod_{i=1}^d\left\lceil\frac{q^{R+1}}{\lambda_i(\Lambda)}\right\rceil.
\end{align*}
In particular, if $D_{\Lambda}\subseteq B(0,q^R)$, then, $\left|\mathcal{D}_{\Lambda}\cap S\cap B(0,q^R)\right|=\left|S\cap \mathcal{D}_{\Lambda}\right|$, so that
\begin{equation}
    \left|B(0,q^R)\cap S\right|=q^N\prod_{i=1}^d\left\lceil \frac{q^{R+1}}{\lambda_{i,\mathcal{C}}(\Lambda)}\right\rceil.
\end{equation}
\end{proof}

\begin{remark}\label{rem:ddddzz}
Let $\mathcal{C}$ be a convex body and let $S$ be a periodic lattice such that $S+\mathcal{C}$ is lattice packing. For $R$ large enough, such that $\mathcal{C}+\mathcal{D}_{\Lambda}\cap S \subset B(0,q^R)$, we have
\begin{equation}\label{eqn:B(0,R)capPack}
        B(0,q^R)\cap (\mathcal{C}+ S)=\mathcal{C}+S\cap B(0,q^R).
    \end{equation}
\end{remark}
As a consequence, we obtain Proposition \ref{prop:Mink}
\begin{proof}[Proof of Proposition \ref{prop:Mink}]
    Let $S$ be a $(\Lambda,q^N)$ periodic lattice, and let $\mathcal{C}$ be a convex body with $m(\mathcal{C}+\mathcal{D}_{\Lambda}\cap S)>\frac{\det(\Lambda)}{q^{d}}$. Then, $S+\mathcal{C}$ is a periodic lattice packing, and therefore, $\delta(S+\mathcal{C})\leq 1$. Since $(\mathbf{s}_1+\mathcal{C})\cap (\mathbf{s}_2+\mathcal{C})=\emptyset$ for distinct $\mathbf{s}_1,\mathbf{s}_2\in S$, by Lemma \ref{lem:PerLattCapBallCnt}, for every $R>\max\left\{\lambda_{d,\mathcal{C}}(\Lambda),\lambda_i(\Lambda)\right\}$, we have
    \begin{equation}
        m\left(B(0,q^R)\cap (S+\mathcal{C})\right)=m(\mathcal{C})\cdot \left|B(0,q^R)\cap S\right|=q^{N+d(R+1)}\frac{m(\mathcal{C})}{\det(\Lambda)}.
    \end{equation}
    Thus
    \begin{equation*}
        1\geq \frac{m\left(B(0,q^R)\cap (S+\mathcal{C})\right)}{q^{Rd}}\rightarrow \frac{q^{N+d}}{\det(\Lambda)}m\left(\mathcal{C}\right),
    \end{equation*}
    as $R \rightarrow \infty$.
    Therefore, $m\left(\mathcal{C}\right)\leq\frac{\det(\Lambda)}{q^{N+d}}$, a contradiction. Hence, $\mathcal{C}$ must contain a non-zero point of $S$.

\end{proof}
\begin{proof}[Proof of Proposition \ref{lem:PackDens}]
    Let $\mathcal{C}'=x^{-1+\log_q\lambda_{1,\mathcal{C}}(S)}\mathcal{C}$. Then, by Lemma \ref{lem:PackRadComp}, $S+\mathcal{C}'$ is a periodic lattice packing. Therefore, by Remark \ref{rem:ddddzz}, for every $R$ large enough, we have
    \begin{equation*}
        B(0,q^R)\cap (\mathcal{C}'+ S)=\mathcal{C}'+B(0,q^R)\cap S.
    \end{equation*}
     Thus, as $R \rightarrow \infty$, Theorem \ref{thm:Mink2nd} and Lemma \ref{lem:PerLattCapBallCnt} imply that
    \begin{align*}
        m\left(\left(S+\mathcal{C}'\right) \cap B(0,q^R)\right)&=|B(0,q^R)\cap S| \cdot  m(\mathcal{C}')\\
        &=\frac{q^{N+(R+1)d}}{\det(\Lambda)}q^{-d}\lambda_{1,\mathcal{C}}(S)^dm(\mathcal{C})\\
        &=q^{N+Rd}\lambda_{1,\mathcal{C}}(S)^d\frac{m(\mathcal{C})}{\det(\Lambda)}.
    \end{align*}
    In conclusion,
    \begin{equation}
        \delta(\mathcal{C}'+S)=\limsup_{R\rightarrow \infty}\frac{m\left((S+\mathcal{C}')\cap B(0,q^R)\right)}{m\left(B(0,q^R)\right)}=q^{N}\lambda_{1,\mathcal{C}}(S)^d\frac{m(\mathcal{C})}{\det(\Lambda)}.
    \end{equation}
\end{proof}
\subsection{Bounds on the Successive Minima}
\begin{proof}[Proof of Theorem \ref{thm:SuccMinBNDS}(\ref{thm:prodBNDPer})]
    Since $\lambda_{i,\mathcal{C}}(\Lambda(\alpha,q^N))=\lambda_i(h^{-1}\Lambda(h^{-1}\alpha,q^N))$, where $\mathcal{C}=h\mathcal{O}^d$, it suffices to prove Theorem \ref{thm:SuccMinBNDS}(\ref{thm:prodBNDPer}) for $\mathcal{C}=\mathcal{O}^d$. Hence, we now assume that $\mathcal{C}=\mathcal{O}^d$.
    
    Let $\mathbf{w}_1,\dots,\mathbf{w}_d\in S$ be linearly independent vectors satisfying $\lambda_i(S)=\Vert \mathbf{w}_i\Vert$. Let $\mathcal{D}=\mathfrak{m}\mathbf{w}_1+\dots+\mathfrak{m}\mathbf{w}_d$. First, we prove that $S+\mathcal{D}$ is a periodic lattice packing. Since $S$ is additive, it suffices to prove that $S \cap \mathcal{D}=\{0\}$. For that, let ${\bf s} \in S \cap \mathcal{D}$ and write 
    \begin{equation}
        \mathbf{s}=\sum_{i=1}^d\gamma_i\mathbf{w}_i,
    \end{equation} where $\gamma_1,\dots,\gamma_d\in \mathfrak{m}$.
    We first prove that for every $i=1,\dots,d$, we have $\Vert \mathbf{s}\Vert\leq q^{-1}\lambda_i(\Lambda)$. This implies that $\mathbf{s}=0$. Firstly, for $i=d$, Lemma \ref{lem:UM} implies that
    \begin{equation}
        \Vert \mathbf{s}\Vert\leq \max_{i=1,\dots,d}\vert \gamma_i\vert\lambda_i(S)\leq q^{-1}\lambda_d(S).
    \end{equation}
    Hence, $\gamma_d=0$. Now, assume that for every $j=i+1,\dots,d$, we have $\gamma_j=0$, and we prove that $\gamma_i=0$. By Lemma \ref{lem:UM}, we have
    \begin{equation}
        \Vert \mathbf{s}\Vert=\left\Vert\sum_{j=1}^i\gamma_j\mathbf{w}_j\right\Vert\leq \max_{j=1,\dots,i}\vert \gamma_j\vert\lambda_j(\Lambda)\leq q^{-1}\lambda_i(\Lambda).
    \end{equation}
    As a consequence, $\gamma_i=0$. By continuing this induction, we obtain that $\mathbf{s}=0$. Therefore $S+\mathcal{D}$ is a periodic lattice packing. If $R>\lambda_d(\Lambda)$, then, by Remark \ref{rem:ddddzz}, we have
    \begin{equation}
        B(0,R)\cap (\mathcal{D}+S)=\mathcal{D}+B(0,R)\cap S.
    \end{equation}
    By Lemma \ref{lem:PerLattCapBallCnt}, we obtain
    \begin{equation}
        \delta(S+\mathcal{D})=\limsup_{R\rightarrow \infty}\frac{m(\mathcal{D})q^{N+(R+1)d}}{q^{R d}\det(\Lambda)}=q^N\frac{\prod_{i=1}^d\lambda_i(S)}{\det(\Lambda)}
    \end{equation}
    Since $\delta(S+\mathcal{D})\leq 1$, if follows that 
    \begin{equation}
        \prod_{i=1}^d\lambda_i(S)\leq \frac{\det(\Lambda)}{q^N}.
    \end{equation}
\end{proof}
\begin{proof}[Proof of Theorem \ref{thm:SuccMinProd}]
Firstly, it suffices to assume that $\mathcal{C}=\mathcal{O}^d$, since $\lambda_{i,\mathcal{C}}(\Lambda(\alpha,q^N))=\lambda_i(h^{-1}\Lambda(h^{-1}\alpha,q^N))$ and $d(\Lambda(\alpha,q^N))=d(h^{-1}\Lambda(h^{-1}\alpha,q^N))$.

    The right inequality of \eqref{eqn:succMinIneq} is a consequence of Theorem \ref{thm:SuccMinBNDS}(\ref{thm:prodBNDPer}) and Lemma \ref{lem:intersectdlamb}. Therefore it suffices to prove the left inequality of \eqref{eqn:succMinIneq}. Let $\mathbf{v}_1,\dots,\mathbf{v}_d$ be a basis of successive minima for $\Lambda$ and let $\mathbf{v}\in \Lambda$. Firstly, write $\alpha=\sum_{i=1}^d\alpha_i\mathbf{v}_i$ and $\mathbf{v}=\sum_{i=1}^da_i\mathbf{v}_i$, where $\alpha_i\in \mathcal{K}_{\infty}$ and $a_i\in \mathcal{R}$. Then,
\begin{equation}
\begin{split}
    \Vert Q\alpha-\mathbf{v}\Vert=\left\Vert\sum_{i=1}^d(Q\alpha_i-a_i)\mathbf{v}^{(i)}\right\Vert=\max_{i=1,\dots,d}\vert Q\alpha_i-a_i\vert\cdot\Vert\mathbf{v}^{(i)}\Vert\\
    \geq\max_{i=1,\dots,d}\vert \langle Q\alpha_i\rangle\vert \lambda_i(\Lambda)=\left\Vert\sum_{i=1}^d\langle Q\alpha_i\rangle\mathbf{v}^{(i)}\right\Vert.
    \end{split}
\end{equation}
As a consequence, 
\begin{equation}
    \min_{\mathbf{v}\in \Lambda}\Vert Q\alpha-\mathbf{v}\Vert=\left\Vert\sum_{i=1}^d\langle Q\alpha_i\rangle\mathbf{v}^{(i)}\right\Vert.
\end{equation}
Assume that $\mathcal{D}_{\Lambda}\cap \Lambda(\alpha,q^N)$ contains $k$ linearly independent vectors, where $1\leq k\leq d$. Then, there exist polynomials $Q_1,\dots,Q_k$ with $\deg(Q_j)\leq N$, such that the vectors $\left\{\sum_{i=1}^d\langle Q_j\alpha_i\rangle\mathbf{v}^{(i)}:j=1,\dots,k\right\}$ are linearly independent. Moreover, there exists a set of indices $I=\{i_1<\dots<i_{d-k}\}\subseteq \{1,\dots,d\}$ for which the set of vectors 
\begin{equation}\label{eqn:IndVects}\mathbf{v}^{(i_1)},\dots,\mathbf{v}^{(i_{d-k})},\sum_{i=1}^d\langle Q_1\alpha_i\rangle\mathbf{v}^{(i)},\dots,\sum_{i=1}^d\langle Q_k\alpha_i\rangle\mathbf{v}^{(i)}\end{equation}
are linearly independent. By the definition of successive minima, there exists $I=\{1\leq i_1<\dots<i_{d-k}\leq d\}$ such that the set \eqref{eqn:IndVects} forms a set of successive minima for $\Lambda(\alpha,q^N)$. Let $\{\ell_1<\dots<\ell_k\}$ be an ordered set of indices such that $\{1,\dots,d\}=I\cup \{\ell_1,\dots,\ell_k\}$. Thus, by \cite[Lemma 2.4]{RW} and Theorem \ref{thm:Mink2nd},
\begin{equation}
\begin{split}
    \prod_{i=1}^d\lambda_i(\Lambda(\alpha,q^N))=\left\Vert\mathbf{v}^{(i_1)}\wedge\dots\wedge\mathbf{v}^{(i_{d-k})}\wedge\sum_{i=1}^d\langle Q_1\alpha_i\rangle\mathbf{v}^{(i)}\wedge\dots\wedge\sum_{i=1}^d\langle Q_k\alpha_i\rangle\mathbf{v}^{(i)}\right\Vert\\
    =\left\Vert\mathbf{v}^{(i_1)}\wedge \dots\wedge\mathbf{v}^{(i_{d-k})}\wedge\sum_{s=1}^k\langle Q_1\alpha_{\ell_s}\rangle \mathbf{v}^{(\ell_s)}\wedge\dots\wedge\sum_{s=1}^k\langle Q_k\alpha_{\ell_s}\rangle\mathbf{v}^{(\ell_s)}\right\Vert\\
    =\Vert \mathbf{v}^{(1)}\wedge\dots\wedge\mathbf{v}^{(d)}\Vert\cdot\left|\det\begin{pmatrix}
        \langle Q_{1}\alpha_{\ell_1}\rangle&\dots&\langle Q_k\alpha_{\ell_1}\rangle\\
        \vdots&\dots&\vdots\\
        \langle Q_{1}\alpha_{\ell_k}\rangle&\dots&\langle Q_k\alpha_{\ell_k}\rangle
    \end{pmatrix}\right|\\
    \geq d(\Lambda(\alpha,q^N))\prod_{i=1}^d\lambda_i(\Lambda)=\det(\Lambda)d(\Lambda(\alpha,q^N)).
\end{split}
\end{equation}
\end{proof}
\subsection{Covering Radii}
We first state the following obvious fact, which is useful in this paper. 
\begin{lemma}
\label{lem:CovRadleq}
     Let $\Lambda$ be a lattice, let $\mathbf{v}_1,\dots,\mathbf{v}_d$ be a basis of successive minima for $\Lambda$, let $N\in \mathbb{N}$, and let $\alpha\in \mathcal{K}_{\infty}^d$. Then, $\operatorname{CovRad}(\Lambda(\alpha,q^N))\leq \operatorname{CovRad}(\Lambda)$.
\end{lemma}
\begin{lemma}
\label{lem:CovRadFundDom}
    Let $S$ be a $(\Lambda,q^N)$ periodic lattice. Let $r>0$ be the covering radius of $S\cap \mathcal{D}_{\Lambda}$ within $\mathcal{D}_{\Lambda}$, that is $r$ is the smallest number such that $\mathcal{D}_{\Lambda}\subseteq S\cap\mathcal{D}_{\Lambda}+B(0,r)$. Then, $\operatorname{CovRad}(S)=r$.
\end{lemma}
\begin{proof}
    Since $\mathcal{D}_{\Lambda}\subseteq B(0,q^{-1}\lambda_d(\Lambda))$, then $r\leq q^{-1}\lambda_d(\Lambda)$. Let $\mathbf{v}\in \mathcal{K}_{\infty}^d$. Then, there exists $\mathbf{u}\in \Lambda$, such that
    \begin{equation}
        \mathbf{v}\in \mathcal{D}_{\Lambda}+\mathbf{u}\subseteq \mathcal{D}_{\Lambda}\cap S+B(0,r)+\mathbf{u}\subseteq \mathcal{D}_{\Lambda}\cap S+\Lambda+B(0,r)=S+B(0,r).
    \end{equation}
    Hence, $\operatorname{CovRad}(S)\leq r$. On the other hand, there exists $\mathbf{v}\in D_{\Lambda}$ for which $\mathbf{v}\notin S\cap \mathcal{D}_{\Lambda}+B(0,q^{-1}r)$. Since $\mathbf{v}$ does not belong to any translate of $\mathcal{D}_{\Lambda}$ by a non-zero $\mathbf{u}\in \Lambda$, then, $\mathbf{v}\notin S+B(0,q^{-1}r)$. Hence, by Lemma \ref{lem:CovRadleq}, $\operatorname{CovRad}(S)=r$.
\end{proof}
\begin{proof}[Proof of Proposition \ref{prop:CovRad}]
    Since $\operatorname{CovRad}_{\mathcal{C}}(\Lambda(\alpha,q^N))=\operatorname{CovRad}(h^{-1}\Lambda(h^{-1}\alpha,q^N))$, where $\mathcal{C}=h\mathcal{O}^d$, then, it suffices to prove Proposition \ref{prop:CovRad} for $\mathcal{C}=\mathcal{O}^d$. By Lemma \ref{lem:CovRadFundDom}, it suffices to compute the covering radius of $\mathcal{D}_{\Lambda}\cap \Lambda(\alpha,q^N)$ in $\mathcal{D}_{\Lambda}$. If $\mathbf{u}\in \mathcal{D}_{\Lambda}$, then, there exist $u_1,\dots,u_d\in \mathfrak{m}$, such that $\mathbf{u}=\sum_{i=1}^du_i\mathbf{v}^{(i)}$. Thus, by Theorem \ref{thm:Mink2nd}, for every $a_1,\dots,a_d\in \mathcal{R}$, we have
    \begin{equation}
    \begin{split}
        \left\Vert\sum_{i=1}^d(u_i-Q\alpha_i-a_i)\mathbf{v}^{(i)}\right\Vert=\max_{i=1,\dots,d}\left|u_i-Q\alpha_i-a_i\right|\lambda_i(\Lambda) \geq \max_{i=1,\dots,d}\left|\langle u_i-Q\alpha_i\rangle\right|\lambda_i(\Lambda)
    \end{split}
    \end{equation}
    Therefore, 
    \begin{equation}
        \min_{\deg(Q)\leq N,\mathbf{v}\in \Lambda}\Vert \mathbf{u}-\mathbf{v}-Q\alpha\Vert=\max_{i=1,\dots,d}\vert \langle u_i-Q\alpha_i\rangle\vert\lambda_i(\Lambda),
    \end{equation}
    so that 
    \begin{equation}
        \operatorname{CovRad}(\Lambda(\alpha,q^N))=\max_{\mathbf{u}\in \mathcal{D}_{\Lambda}}\min_{\deg(Q)\leq N,\mathbf{v}\in \Lambda}\Vert \mathbf{u}-\mathbf{v}-Q\alpha\Vert=\max_{\mathbf{u}\in \mathcal{D}_{\Lambda}}\max_{i=1,\dots,d}\vert \langle u_i-Q\alpha_i\rangle\vert\lambda_i(\Lambda).
    \end{equation}
    Hence, $\operatorname{CovRad}(\Lambda(\alpha,q^N))<q^{-\ell}$ if and only if for every $\mathbf{u}=\sum_{i=1}^du_i\mathbf{v}^{(i)}\in \mathcal{D}_{\Lambda}$, there exists $Q=\sum_{j=0}^NQ_jx^j\in \mathcal{R}_{\leq N}$ such that for every $i=1,\dots,d$ such that $\lambda_i(\Lambda)>q^{-\ell}$, we have $\vert\langle u_i-Q\alpha_i\rangle\vert<\lambda_i(\Lambda)^{-1}q^{-\ell}<q^{-(\ell+\log_q\lambda_i(\Lambda))}$. Let $\alpha_i=\sum_{j=1}^{\infty}\alpha_{j,i}x^{-j}$ and $u_i=\sum_{j=1}^{\infty}u_{j,i}x^{-j}$. If $\langle Q\alpha_i\rangle=\sum_{j=1}^{\infty}c_jx^{-j}$, then, 
    \begin{equation}
        c_j=Q_0\alpha_{j,i}+Q_1\alpha_{j+1,i}+\dots+Q_N\alpha_{j+N,i}.
    \end{equation}
    Since $\vert \langle u_i-Q\alpha_i\rangle\vert<q^{-\ell}$ for every $i=1,\dots, d$ satisfying $\ell+\log_q\lambda_i(\Lambda)>0$, then, for every such $i$, we have
    \begin{equation}
        \Delta_{\alpha_i}(\ell+\log_q\lambda_i(\Lambda),N+1)\begin{pmatrix}
            Q_0\\
            \vdots\\
            Q_N
        \end{pmatrix}=\begin{pmatrix}
            u_{1i}\\
            \vdots\\
            u_{\ell+\log_q\lambda_i(\Lambda),i}
        \end{pmatrix}.
    \end{equation}
    Thus, if $\ell+\log_q\lambda_i(\Lambda)>0$, then, $\operatorname{rank}(\Delta_{\alpha_i}(\ell+\log_q\lambda_i(\Lambda),N+1))=\ell+\log_q\lambda_i(\Lambda)$. Let $j$ be the smallest index such that $\lambda_i(\Lambda)>q^{-\ell}$. As a consequence, the matrix $$\Delta_{\alpha_1,\dots,\alpha_d}(\ell+\log_q\lambda_1(\Lambda),\dots,\ell+\log_q\lambda_d(\Lambda);N+1)=\begin{pmatrix}
        \Delta_{\alpha_j}(\ell+\log_q\lambda_j(\Lambda),N+1)\\
        \vdots\\
        \Delta_{\alpha_d}(\ell+\log_q\lambda_d(\Lambda),N+1)
    \end{pmatrix}$$ 
    has rank $\sum_{i=j}^d(\ell+\log_q\lambda_i(\Lambda))=\sum_{i=1}^d\max\{\ell+\log_q\lambda_i(\Lambda),0\}$. Hence, we have
    \begin{small}
    \begin{equation}
        \operatorname{CovRad}(\Lambda(\alpha,q^N))=q^{-\left(1+\max\left\{\ell:\operatorname{rank}(\Delta_{\alpha_1,\dots,\alpha_d}(\ell+\log_q\lambda_1(\Lambda),\dots,\ell+\log_q\lambda_d(\Lambda);N+1))=\sum_{i=1}^d\max\{\ell+\log_q\lambda_i(\Lambda),0\}\right\}\right)}.
    \end{equation}
    \end{small}.
\end{proof}
\section{Acknowledgements}
I would like to thank Angelot Behajaina for introducing me to these questions about periodic lattices and for useful discussions about geometry of numbers. I would also like to thank the anonymous referees for their useful comments which helped improve the presentation of this manuscript and improve some inaccuracies in the earlier version.
\bibliography{Ref}
\bibliographystyle{amsalpha}
\end{document}